\documentclass[12pt]{amsart}
\usepackage{amsmath,amscd,amssymb,amsfonts}
\newcommand{\h}{\hbox}
\newcommand{\q}{\quad}

\newcommand{\bs}{\par\bigskip}
\newcommand{\ms}{\par\medskip}
\newcommand{\sk}{\par\smallskip}
\newcommand{\bsn}{\par\bigskip\noindent}
\newcommand{\msn}{\par\medskip\noindent}
\newcommand{\skn}{\par\smallskip\noindent}
\newcommand{\mprod}{\h{$\prod$}}
\newcommand{\msum}{\h{$\sum$}}
\newcommand{\mopl}{\h{$\bigoplus$}}

\newcommand{\mcup}{\h{$\bigcup$}}
\newcommand{\ssb}{\raise.15ex\h{${\scriptscriptstyle\bullet}$}}
\newcommand{\ssc}{\,\raise.15ex\h{${\scriptstyle\circ}$}\,}
\newcommand{\Ab}{{\mathbb A}}
\newcommand{\C}{{\mathbb C}}
\newcommand{\N}{{\mathbb N}}
\newcommand{\Ob}{{\mathbb O}}
\newcommand{\PS}{{\mathbb P}{\mathbb S}}
\newcommand{\Q}{{\mathbb Q}}

\newcommand{\Z}{{\mathbb Z}}
\newcommand{\sbo}{{\bf s}}
\newcommand{\ob}{{\bf 1}}

\newcommand{\B}{{\mathcal B}}

\newcommand{\E}{{\mathcal E}}
\newcommand{\Hc}{{\mathcal H}}
\newcommand{\G}{{\mathcal G}}
\newcommand{\OO}{{\mathcal O}}
\newcommand{\Ct}{\widetilde{C}}
\newcommand{\Dt}{\widetilde{D}}
\newcommand{\dd}{\partial}
\newcommand{\dt}{\dd_t}
\newcommand{\dti}{\dd_t^{-1}}
\newcommand{\dtip}{\{\!\{\dti\}\!\}}
\newcommand{\dtipdt}{\dtip[\dd_t]}

\newcommand{\ee}{\widetilde{e}}
\newcommand{\gt}{\widetilde{g}}
\newcommand{\hh}{\widetilde{h}}
\newcommand{\sss}{\widetilde{s}}
\newcommand{\SSS}{\widetilde{S}}
\newcommand{\ttt}{\widetilde{t}}
\newcommand{\uu}{\widetilde{u}}
\newcommand{\vv}{\widetilde{v}}
\newcommand{\pit}{\widetilde{\pi}}
\newcommand{\PSh}{\widehat{\PS}}
\newcommand{\Gr}{{\rm Gr}}
\newcommand{\into}{\hookrightarrow}
\newcommand{\onto}{\mathop{\rlap{$\to$}\hskip2pt\h{$\to$}}}
\newcommand{\simto}{\buildrel\sim\over\longrightarrow}
\newcommand{\ges}{\geqslant}
\newcommand{\les}{\leqslant}
\newcommand{\bl}{\bigl}
\newcommand{\br}{\bigr}
\newcommand{\1}{\hskip1pt}
\newcommand{\fa}{\raise1pt\hbox{$\forall$}\1}

\newcommand{\ep}{\varepsilon}
\newcommand{\al}{\alpha}
\newcommand{\nee}{\h{\rlap{$=$}\raise1pt\h{$\scriptscriptstyle \,\,/$}}}
\newcommand{\precnee}{\,\h{$\rlap{$\prec$}\raise-5pt\nee$}\,\,\,}
\newcommand{\succnee}{\,\h{$\rlap{$\succ$}\raise-5pt\nee$}\,\,\,}
\begin{document}
\title{Deformations of abstract Brieskorn lattices}
\author{Morihiko Saito}
\begin{abstract}
We study certain deformations of abstract Brieskorn lattices in fixed abstract Gauss-Manin systems, and show that the ambiguity of expressions of deformations coming from automorphisms of base spaces is essentially the same as the one coming from the choice of opposite filtrations, and hence is finite dimensional, although the freedom of parameters in the expressions of deformations is infinite dimensional. As a consequence, we can prove the non-existence of a versal deformation of the Fourier transform of this abstract Brieskorn lattice with expected dimension. This shows that the generation condition is quite essential for the existence of versal deformations with expected dimensions in the absolute case.
\end{abstract}
\maketitle
\centerline{\bf Introduction}
\bsn
Consider an abstract Brieskorn lattice $\B$ generated by $\dt^{-j}e_j$ ($j=0,\dots,r$) over $\C\dtip$ in an abstract Gauss-Manin system $\G$ (see \cite{bl}), that is,
$$\B:=\msum_{j=0}^r\,\C\dtip\dt^{-j}\1 e_j\,\subset\,\G:=\mopl_{j=0}^r\,\C\dtipdt\1 e_j.
\leqno(1)$$
Here $r$ is an integer at least $2$, and the action of $t$ is given by
$$t\1\dt^{-k}e_j=(k+1)\dt^{-k-1} e_j\,\,\,(k\in\Z)\,\,\,\h{independently of}\,\,\,\,j\in\{0,\dots,r\}.
\leqno(2)$$
\sk
Let $S$ be a sufficiently small open polydisk in $\C^r$ with coordinates $s_1,\dots,s_r$. Let $h(s_1)$ be a holomorphic function of $s_1$ with $h(0)=h'(0)=0$, that is,
$$h(s_1)=\msum_{k\ges 2}\,c_k\1s_1^k\,\,\,(c_k\in\C)\q\h{with}\q\msum_{k\ges 2}\,|c_k|\1\ep^k<\infty\,\,\,\,(\exists\,\ep>0).$$
Here we set $\,h'(s_1):=\dd_{s_1}h(s_1)$, $h''(s_1):=\dd_{s_1}^2h(s_1)$.
\sk
Consider a deformation $\B_{h,S}$ of $\B$ inside $\G$ defined by
$$\B_{h,S}:=\msum_{j=0}^r\,\OO_S\dtip\1v_j\,\subset\,\G_S:=\G\,\widehat{\otimes}_{\C}\,\OO_S=\mopl_{j=0}^r\,\OO_S\dtipdt\1e_j,
\leqno(3)$$
with
$$\dd_{s_i}e_j=0\q(\fa i,j),$$
and
$$\aligned v_0:={}&e_0+\msum_{j=1}^r\,s_j\dt^{1-j}e_j+h(s_1)\1 e_2,\\ v_j:={}&\dti\dd_{s_j}v_0=\begin{cases}\dt^{-1}e_1+h'(s_1)\1 \dt^{-1}e_{2}&\h{if}\,\,\,j=1,\\ \dt^{-j}e_j&\h{if}\,\,\,j\ges 2.\end{cases}\endaligned
\leqno(4)$$
Here $v_0$ is a ``primitive form" of this deformation $\B_{h,S}$ associated with the opposite filtration $U$ defined by the $e_j$, see \cite{SK}, \cite{bl} for the geometric case. (This notion was quite useful for the construction of this example.) By (2) we have
$$\aligned t\1 v_0&=\dti v_0+\msum_{j=2}^r\,(j-1)\1 s_j\1 v_j,\\ t\1 v_j&=(j+1)\1\dti v_j\q\,\,(j=1,\dots,r),\\ \dti\dd_{s_i}v_0&=v_i\1\q\q\q\q\q\q(i=1,\dots,r),\\ \dti\dd_{s_i}v_j&=\begin{cases}h''(s_1)\1 v_2\q\,\,\h{}&\h{if}\,\,\,(i,j)=(1,1),\\ 0&\h{if}\,\,\,j\ges 1,\,\,\,(i,j)\ne(1,1).\end{cases}\endaligned
\leqno(5)$$
So $\B_{h,S}$ is an $\E(0)$-{\it submodule} of the $\E$-module $\G_S$, that is, $\B_{h,S}$ is stable by the actions of $t$, $\dti$, $\dti\dd_{s_j}$ ($j=1,\dots,r$), see \cite{SKK}. We say that $\B_{h,S}$ is a {\it deformation of special type} associated with a holomorphic function $h$ (or more precisely, a convergent power series $h$, since the coordinate $s_1$ is {\it fixed}). Note that $\G_S$ is a {\it constant} deformation of $\G$. (This corresponds to a $\mu$-constant deformation when the Brieskorn lattices are associated with hypersurface isolated singularities, see \cite{Br}, \cite{Ph}, \cite{bl}, etc.)
\sk
The above function $h$ is {\it not} uniquely determined by the deformation $\B_h$, since there is an {\it ambiguity} of $h$ coming from {\it automorphisms} of $S$ and $\G_S$. We study this ambiguity by using the {\it uniqueness} of the generators $v_j$ associated with an opposite filtration $U$ (see Corollary~(1.4) below), and prove the following.
\msn
{\bf Theorem~1.} {\it There is only finite dimensional $($more precisely, $3$-dimensional\1$)$ ambiguity of $h$ under automorphisms of $S$ and $\G_S$. This ambiguity is the same as the one coming from the choice of the opposite filtration $U$, which is given by the $e_j$ in this case.}
\ms
Note that the moduli (or freedom) of $h$ is {\it infinite dimensional,} using the Taylor expansion as above, since there is no condition on $h$ except that $h(0)=h'(0)=0$.
\sk
To show Theorem~1, let $\SSS$ be a sufficiently small open polydisk in $\C^r$ with coordinates $\sss_1,\dots,\sss_r$. Define $\G_{\SSS}$ with $S$ replaced by $\SSS$, where the $e_j$ are denoted by the $\ee_j$ ($j=0,\dots,r$). We have $\ttt=t$, since deformations of $\G$, $\B$ are considered. Let $\hh(\sss_1)$ be a holomorphic function of $\sss_1$. Define $\vv_j$ ($j=0,\dots,r$) and $\B_{\hh,\SSS}\subset\G_{\SSS}$ as above with $h,e_j,s_i$ replaced by $\hh,\ee_j,\sss_i$. Assume there is an isomorphism
$$\rho:(S,0)\simto(\SSS,0),$$
together with an isomorphism of $\OO_S\dtipdt$-modules compatible with the actions of $t,\dd_{s_i}$ and preserving lattices:
$$\phi:\bl(\G_S,\B_{h,S}\br)\simto\rho^*\bl(\G_{\SSS},\B_{\hh,\SSS}\br).
\leqno(6)$$
Restricting to the intersection of the kernels of the actions of $t$, $\dd_{s_i}$ (which coincides with $\sum_{j=0}^r\C\1\dt e_j$), the last isomorphism can be determined by the $\C$-linear map
$$A:H=\mopl_{j=0}^r\,\C\1e_j\simto\mopl_{j=0}^r\,\C\1\ee_j,$$
corresponding to a matrix $(A_{i,j})$ of size $(r+1)\times(r+1)$ with $\C$-coefficients. We will identify the $\ee_j$ with elements of $H$ via the inverse of the linear map $A$ so that
$$e_i=\msum_{j=0}^i\,A_{i,j}\,\ee_j\q\h{in}\,\,\,H.
\leqno(7)$$
Here we have $A_{i,j}=0$ for $i<j$, since $A$ preserves the Hodge filtration $F_0$ on $H$ at the origin of $S,\SSS$ (which is defined by $\Gr_V^1(\dt^{-p}\B)$, see (1.1) below). We assume $A_{i,i}=1$ to avoid problems coming from the ambiguity of the $v_j$ by non-zero constant multiples.
So the matrix $(A_{i,j})$ corresponds to the choice of the opposite filtration $U$ to the Hodge filtration $F_0$ on $H$ at the origin of $S$.
\sk
We have the convergent direct sum decomposition
$$\G=\widehat{\mopl}_{k\in\Z}\,\G^k,\q\G_S=\widehat{\mopl}_{k\in\Z}\,\G^k_S,
\leqno(8)$$
with
$$\G^k:={\rm Ker}(\dt t-k)\subset\G,\q\G^k_S:={\rm Ker}(\dt t-k)\subset\G_S.$$
We have the canonical isomorphisms
$$\G^k_S\buildrel{\dt^{k-1}}\over\longrightarrow\G^1_S=\Hc_S:=H\otimes_{\C}\OO_S=\mopl_j\,\OO_S\1e_j\,\subset\,\G_S\q(k\in\Z),$$
and similarly for $\G^k$ by omitting $S$ and replacing $\OO_S$ with $\C$.
These give the canonical projections
$$pr_S^{(k)}:\G_S=\widehat{\mopl}_{k\in\Z}\,\G^k_S\onto\G^k_S=\Hc_S\q(k\in\Z).$$
We have moreover the canonical projections
$$\pit_j:\Hc_S=H\otimes_{\C}\OO_S\to\OO_S\q(j=0,\dots,r),$$
which are the scalar extension by $\C\into\OO_S$ of the projections $H\to\C\1\ee_j=\C$. Recall that the $\ee_j$ are identified with elements of $H$ via the inverse of the linear map $A$. Let
$$v'_j\in\B_{h,S}\q(j=0,\dots,r)$$
be the generators of $\B_{h,S}$ associated to the basis $(\ee_j)$ (which determines the opposite filtration $U$), see Corollary~(1.4) below. Define
$$g_j:=\pit_j\ssc pr_S^{(j)}(v'_0)\,\,\,(j=1,\dots,r),\q\gt:=\pit_2\ssc pr_S^{(1)}(v'_0)\,\,\in\,\,\OO_S.
\leqno(9)$$
These are uniquely determined by $\B_{h,S}$ and the $\ee_j$ (or the matrix $(A_{i,j})$).
\sk
By the uniqueness of generators associated with the basis $(\ee_j)$ (see Corollary~(1.4) below), we get the following.
\msn
{\bf Theorem~2.} {\it There are equalities in} $\,\OO_S:$
$$g_j=\rho^*\sss_j\,\,\,(j=1,\dots,r),\q\gt=\rho^*\hh.$$
\ms
As a corollary, the morphism $\rho$ is uniquely determined by the $g_j$. Theorem~1 follows from Theorem~2 except for the dimension of the ambiguity (see (2.2) below for the latter).
\sk
Theorem~1 together with a remark after it implies a certain difficulty in constructing a versal deformation of this Brieskorn lattice $\B$ {\it inside\1} the fixed Gauss-Manin system $\G$. For instance, its base space $S$ cannot be $r$-dimensional. Here one cannot simply replace $h(s_1)$ in the definition of $v_0$ by a {\it new\1} coordinate $s_{r+1}$, since one must have an $\E(0)$-submodule which is stable by the action of the $\dti\dd_{s_i}$ ($j=1,\dots,r+1$) and $\dim\Gr_1^{F_0}H=1$. As a consequence, we can prove the non-existence of a versal deformation of the Fourier transform of $\B$ whose base space has dimension $r+1$, see (3.8) below. This shows that the generation condition (see \cite{HM}, \cite{Sab}) is quite essential for the existence of versal deformations, see Remark~(2.5) below. It seems quite difficult to extend the argument in this paper to the case where the generation condition is satisfied, see Remarks~(2.3) and (2.4) below.
\sk
This paper is written to answer a question of C.~Hertling. I thank him for this interesting question. This work is partially supported by JSPS Kakenhi 15K04816.
\sk
In Section~1 we generalize the theory of canonical generators associated with an opposite filtration in \cite{bl} to the deformation case.
In Section~2 we prove the main theorems.
In Section~3 we show the non-existence of a versal deformation of the Fourier transform $\B^F$ of $\B$ with dimension $r+1$ using the argument in the proof of Theorem~1.
\bs\bs
\vbox{\centerline{\bf 1. Canonical generators of deformations of Brieskorn lattices}
\bsn
In this section we generalize the theory of canonical generators associated with an opposite filtration in \cite{bl} to the deformation case.}
\msn
{\bf 1.1.~Deformation of abstract Brieskorn lattices.} Let $(S,0)$ be a germ of complex manifold with local coordinates $s_i$. Let $\B_S\subset\G_S$ be a deformation of an abstract Brieskorn lattice $\B$ in a fixed abstract Gauss-Manin system $\G$. Here
$\B_S$, $\G_S$, $\B$, $\G$ are freely generated over $\OO_S\dtip$, $\OO_S\dtip[\dt]$, $\C\dtip$, $\C\dtip[\dt]$ with
$$\OO_S\dtip:=\bl\{\msum_{k\ges 0}\,g_k\dt^{-k}\,\big|\,\msum_{k\ges 0}\,g_ks_0^k/k!\in\OO_{S\times\C}\,\,(g_k\in\OO_S)\br\},
\leqno(1.1.1)$$
and $\B_S$ is stable by the action of $t$, $\dti$, $\dti\dd_{s_i}$.
\sk
We assume that $\G_S$ is a {\it constant\1} deformation of $\G$ so that there are converging direct sum decompositions (related to (1.1.1) as in (1.1.3) below)
$$\G=\widehat{\mopl}_{\al\in\C}\,\G^{\al},\q\G_S=\widehat{\mopl}_{\al\in\C}\,\G^{\al}_S\q\h{with}\q\G^{\al}_S=\OO_S\otimes_{\C}\G^{\al}\,\,\,(\al\in\C),
\leqno(1.1.2)$$
where we set for any integer $m$ greater than or equal to the rank of $\B$ over $\C\dtip$
$$\G^{\al}:={\rm Ker}(\dt t-\al)^m\subset\G,\q\G^{\al}_S:={\rm Ker}(\dt t-\al)^m\subset\G_S.$$
In particular, the $\G^{\al}_S$ are finite free $\OO_S$-modules.
\sk
\vbox{We choose and fix a total order $\prec$ of $\C$ such that we have for any $\al,\beta\in\C$
\skn
(i) $\al\prec\al+1$,
\skn
(ii) $\al\prec\beta\iff\al+1\prec\beta+1$,
\skn
(iii) $\al\prec\beta+m\,$ for some integer $m\gg 0$ (depending on $\al,\beta$).}
\sk
If the monodromy of $\G$ is quasi-unipotent so that the converging direct sum decompositions in (1.1.2) are indexed by $\Q$, we {\it assume\1} that the restriction of $\prec$ to $\Q$ is the natural ordering. Set
$$\Lambda:=\{\al\in\C\mid 0\precnee\al\prec 1\},$$
$$\Hc_S:=\mopl_{\al\in\Lambda}\,\Hc_S^{\al}\q\h{with}\q\Hc_S^{\al}:=\G^{\al}_S\,\,\,\,(\al\in\Lambda).$$
Then
$$\G_S=\OO_S\dtip[\dt]\otimes_{\OO_S}\Hc_S.
\leqno(1.1.3)$$
Note that $\Hc_S^{\al}=0$ except for a finite number of $\al\in\Lambda$.
\sk
Let $V$ be the filtration of Kashiwara and Malgrange indexed decreasingly by $\C$ using the total order $\prec$ so that
$$V^{\al}\G_S:=\widehat{\mopl}_{\beta\succ\al}\,\G^{\beta}_S,\q\q\Gr_V^{\al}\G_S=\G^{\al}_S\q(\al\in\C).
\leqno(1.1.4)$$
We have the {\it Hodge filtration} $F$ on $\Hc_S^{\al}$ defined by
$$F_p\Hc_S^{\al}:=\Gr_V^{\al}\dt^p\B_S\subset\Gr_V^{\al}\G_S=\Hc_S^{\al}\q(\al\in\Lambda).$$
This is a finite filtration, and we have for any $\al\in\C$
$$\dt^{-p}\B_S\subset V^{\al}\G_S\subset\dt^p\B_S\q\h{if}\,\,\,\,p\gg 0\,\,\,\h{(depending on}\,\,\,\al).
\leqno(1.1.5)$$
Indeed, $\B_S$, $V^{\al}\G_S$ are finitely generated over $\OO_S\dtip$, and
$$\mcup_{p\ges 0}\,\dt^p\B_S=\mcup_{\al\in\C}\,V^{\al}\G_S=\G_S.$$
\sk
We {\it assume} the following condition:
$$\h{\it The $\,\Gr_p^F\Hc_S^{\al}\,$ are locally free $\OO_S$-modules $(\forall\,p\in\Z,\,\al\in\Lambda)$.}
\leqno(1.1.6)$$
This is satisfied in the geometric case, since we assume that the restriction of the total order $\prec$ to $\Q$ is the natural order, see \cite{Va2}. (However, it does not necessarily hold without this assumption.)
\sk
We denote by $F_s$ the restriction of the Hodge filtration $F$ to
$$\Hc_{S,s}^{\al}\otimes_{\OO_{S,s}}\C=H^{\al}\,\,\,\h{for}\,\,\, s\in S.$$
This is a finite filtration on finite dimensional vector spaces $H^{\al}$. By (1.1.6), the dimension of $F_{\!s,p}H^{\al}$ is independent of $s\in S$ for any $p\in\Z$, $\al\in\Lambda$.
\msn
{\bf 1.2.~Opposite filtrations.} In the notation and assumption of (1.1), we say that a finite decreasing filtration $U$ on $H^{\al}$ is {\it opposite} to $F_s$ if $U$ gives the direct sum decomposition
$$\mopl_{p\in\Z}\,U^pF_{\!s,p}H^{\al}\simto H^{\al},
\leqno(1.2.1)$$
inducing a {\it splitting} of $F_s$ and $U$ on $H^{\al}$. This is equivalent to each of the following equivalent conditions (see \cite{th2}):
$$\Gr^p_U\Gr_q^{F_s}H^{\al}=0\q\h{if}\q p\ne q.
\leqno(1.2.2)$$
$$F_{\!s,p}\,H^{\al}\oplus U^{p+1}H^{\al}\simto H^{\al}\q(\forall\,p\in\Z).
\leqno(1.2.3)$$
\sk
If $U$ is opposite to $F_0$, then it is opposite to $F_s$ for any $s\in S$ by (1.2.3) (shrinking $S$ if necessary). In this case we say that $U$ is opposite to $F$ on $S$. It defines a filtration $U$ on $\Hc_S^{\al}$ via the scalar extension as well as a filtration $U$ on $\G_S^{\beta}$ for any $\beta\in\C$ so that
$$\dt^k:U^p\G_S^{\beta}\simto U^{p+k}\G_S^{\al}:=U^{p+k}\Hc_S^{\al}\q\h{if}\q\al:=\beta-k\in\Lambda.
\leqno(1.2.4)$$
The Hodge filtration $F$ is also defined on the $\G_S^{\beta}$, and we have the direct sum decompositions
$$\mopl_{p\in\Z}\,U^pF_p\,\G_S^{\beta}\simto\G_S^{\beta}\q(\forall\,\beta\in\C).
\leqno(1.2.5)$$
$$F_p\,\G_S^{\beta}\oplus U^{p+1}\G_S^{\beta}\simto\G_S^{\beta}\q(\forall\,p\in\Z,\,\beta\in\C).
\leqno(1.2.6)$$
\sk
Define the opposite filtration $U$ on $\G_S$ by
$$U^p\G_S:=\mopl_{\beta\in\C}\,U^p\G_S^{\beta}\subset\G_S.$$
Then we have the direct sum decomposition
$$\dt^p\B_S\oplus U^{p+1}\G_S\simto\G_S\q(\forall\,p\in\Z).
\leqno(1.2.7)$$
Indeed, (1.2.6) implies the canonical isomorphism
$$\dt^p\B_S\simto\G_S/U^{p+1}\G_S\q(\forall\,p\in\Z),$$
by taking the graded pieces of the filtration $V$ on it (and using (1.1.5)).
\sk
We have moreover the canonical isomorphism
$$\dt^p\B_S\cap U^p\G_S\simto\Gr_U^p\G_S\q(\forall\,p\in\Z).
\leqno(1.2.8)$$
Indeed, (1.2.7) implies the direct sum decomposition
$$(\dt^p\B_S\cap U^p\G_S)\oplus U^{p+1}\G_S\simto U^p\G_S\q(\forall\,p\in\Z),
\leqno(1.2.9)$$
by using the inclusion $U^{p+1}\G_S\subset U^p\G_S$.
\sk
Generalizing \cite[Proposition 3.4]{bl}, we have the following.
\msn
{\bf Theorem~1.3.} {\it In the notation and assumptions of $(1.1)$, let $U$ be a filtration on $H^{\al}$ opposite to $F$ for any $\al\in\Lambda$. Then, for any $e\in U^pH^{\al}$ with $p\in\Z$, $\al\in\Lambda$, there is a unique $w\in\B_S$, shrinking $S$ if necessary, such that
$$w-\dt^{-p}e\,\in\,\mopl_{\beta\succ\al+p}\,U^1\G_S^{\beta},
\leqno(1.3.1)$$
where the canonical inclusion $H^{\al}\into\Hc_S^{\al}\into\G_S$ is used.}
\msn
{\it Proof.} This follows from (1.2.9) for $p=0$, except for the assertion that the direct sum is indexed by $\beta\succ\al+p$. Indeed, (1.2.8) implies that for $\dt^{-p}e\in U^0\G_S$, there are unique
$$w\in\B_S\cap U^0\G_S,\q u\in U^1\G_S,$$
satisfying
$$w-u=\dt^{-p}e.
\leqno(1.3.2)$$
This implies (1.3.1) except that the direct sum in it is indexed by $\beta\succ\al+p$. For the proof of the last assertion, we have to show
$$u\in V^{\al+p}\G_S.
\leqno(1.3.3)$$
\sk
Assume
$$u\in V^{\beta}\G_S\q\h{for some}\,\,\,\beta\precnee\al+p.$$
By (1.2.6) for $p=1$, we have
$$\Gr_V^{\beta}u=0\q\h{in}\,\,\,\,\Gr_V^{\beta}\G_S.$$
(Indeed, $\Gr_V^{\beta}u\in F_0\G_S^{\beta}\cap U^1\G_S^{\beta}$, since $\Gr_V^{\beta}\dt^{-p}e=0$ in $\Gr_V^{\beta}\G_S$.) So the assertion (1.3.3) follows. This completes the proof of Theorem~(1.3).
\msn
{\bf Corollary~1.4.} {\it In the notation and assumptions of $(1.1)$, let $U$ be a filtration on $H^{\al}$ opposite to $F$ for any $\al\in\Lambda$. Let $(e_{\al,p,j})$ be a basis of $U^pF_{\!0,p}H^{\al}$ for any $p\in\Z$, $\al\in\Lambda$. Then, shrinking $S$ if necessary, there are unique $w_{\al,p,j}\in\B_S$ satisfying
$$w_{\al,p,j}-\dt^{-p}e_{\al,p,j}\,\in\,\mopl_{\beta\succ\al+p}\,U^1\G_S^{\beta}\q\q(\forall\,\al,p,j),
\leqno(1.4.1)$$
and generating $\B_S$ over $\OO_S\dtip$.}
\msn
{\it Proof.} By Theorem~(1.3), it remains to show the last assertion. By (1.1.5) it is enough to show that they generate $\Gr_V^{\ssb}\B_S$ over $\OO_S[\dti]$ after passing to the graded quotients of $V$. Here we may restrict to $s=0$ and pass to $\Gr_V^{\ssb}\B$ using Nakayama's lemma. Then the assertion follows from the assumption that $(e_{\al,p,j})$ is a basis of $U^pF_{\!0,p}H^{\al}$. This completes the proof of Corollary~(1.4).
\msn
{\bf 1.5.~Inductive proof of Theorem~(1.3).} Shrinking $S$ if necessary, we can show by increasing induction on $\gamma\succnee\al+k$ that there are $w_{\gamma}\in\B_S$ together with {\it unique} $u_{\beta}\in U^1\G_S^{\beta}$ for $\al+k\prec\beta\precnee\gamma$ such that
$$w_{\gamma}=\dt^{-k}e+\msum_{\beta}\,u_{\beta}\mod V^{\gamma}\G_S,
\leqno(1.5.1)$$
where the $u_{\beta}$ are independent of $\gamma$ by its uniqueness.
\sk
Indeed, assume (1.5.1) holds for some $\gamma\succnee\al+k$. We may assume $\Gr_V^{\gamma}\G_S\ne 0$. There are $\gamma'\succnee\gamma$, $\gamma''\precnee\gamma$ such that $\Gr_V^{\gamma'}\G_S\ne 0$, $\Gr_V^{\gamma''}\G_S\ne 0$, and
$$\Gr_V^{\gamma}\G_S=V^{\gamma}\G_S/V^{\gamma'}\G_S,\q\Gr_V^{\gamma''}\G_S=V^{\gamma''}\G_S/V^{\gamma}\G_S.$$
By (1.1.12) for $p=0$ together with the isomorphism $\Gr_V^{\gamma}\B_S=F_0\G^{\gamma}_S$, there are
$$u_{\gamma}\in U^1\G_S^{\gamma},\q w'_{\gamma}\in V^{\gamma}\B_S,$$
such that
$$\bl(w_{\gamma}-\dt^{-k}e-\msum_{\beta}\,u_{\beta}\br)-u_{\gamma}-w'_{\gamma}\in V^{\gamma'}\G_S,
\leqno(1.5.2)$$
using the isomorphism (1.2.6) for $p=0$. Then (1.5.1) for $\gamma'$ holds by setting
$$w_{\gamma'}:=w_{\gamma}-w'_{\gamma}\in\B_S.$$
\sk
At the beginning of induction, we have to treat the case where $\gamma$ is sufficiently close to $\al+k$ so that
$$\Gr_V^{\al+k}\G_S=V^{\al+k}\G_S/V^{\gamma}\G_S.$$
Here the assertion (1.5.1) is equivalent to the existence of $w\in F_0U^0\G_S^{\al+k}$ corresponding to $\dt^{-k}e$ by the canonical isomorphisms
$$F_0U^0\G_S^{\al+k}\simto\Gr_U^0\G_S^{\al+k}\,\bl(=\Gr_U^k\Hc_S^{\al}\br).$$
and this follows from (1.2.6) using the same argument as in the proof of (1.2.9).
\sk
These arguments, giving an explicit construction of the $u_{\beta}$, are more useful in some cases, for instance, in the proof of Theorem~1 in (2.2) below.
\bs\bs
\vbox{\centerline{\bf 2. Proof of main theorems}
\bsn
In this section we prove the main theorems.}
\msn
{\bf 2.1.~Proof of Theorem~2.} This follows from Corollary~(1.4).
\ms
(It may be better to use Theorem~(1.3) for the proof of Theorem~2, since only $v'_0$ is used in the definition of $g_j$, $\gt$.)
\msn
{\bf 2.2.~Proof of Theorem~1.} By Theorem~2 together with the arguments in the introduction, it remains to show that the ambiguity of $h$ as a power series of $s_1$ is three dimensional, and that any matrices $(A_{i,j})$ with $A_{i,i}=1$ ($i=0,\dots,r$) and $A_{i,j}=0$ ($i<j$) are allowed here.
\sk
We first show that
$$\dt^{-k}\ee_j\in\B_{h,S}\q\h{if}\q k\ges 2,\,\,\,0\les j\les r,\,\,\,k\ges j.
\leqno(2.2.1)$$
Here we may replace $\ee_j$ with $e_j$ for any $j\in\{0,\dots,r\}$ by (7). Indeed, we have
$$F_{0,k}H=\msum_{j=0}^k\,\C\1 e_j=\msum_{j=0}^k\,\C\1\ee_j\q\h{for any}\,\,\,0\les k\les r.$$
By the definition of $v_j$ ($j=2,\dots,r$), the assertion (2.2.1) with $\ee_j$ replaced by $e_j$ for any $j\in\{1,\dots,r\}$ holds when
$$k\ges 2,\,\,\,2\les j\les r,\,\,\,k\ges j.$$
This implies that
$$\dt^{-k}\bl(e_0+s_1e_1+h(s_1)e_2\br)\in\B_{h,S}\q(k\ges 1).
\leqno(2.2.2)$$
Indeed, (2.2.2) for $k=1$ is proved by considering $\dti v_1$ and using the above assertion. By the definition of $v_1$ we also get
$$\dt^{-k}\bl(e_1+h'(s_1)e_2\br)\in\B_{h,S}\q(k\ges 1).
\leqno(2.2.3)$$
So (2.2.1) follows from these.
\sk
Set now
$$\al:=A_{1,0},\q\beta:=A_{2,0},\q\gamma:=A_{2,1}.$$
By (7) we have
$$e_0=\ee_0,\q e_1=\al\1\ee_0+\ee_1,\q e_2=\beta\1\ee_0+\gamma\1\ee_1+\ee_2.$$
Hence
$$v_0=\bl(1+\al\1 s_1+\beta\1 h(s_1)\br)\1\ee_0+\bl(s_1+\gamma\1 h(s_1)\br)\1\ee_1+h(s_1)\1\ee_2\mod\, V^{>1}\G_S.
\leqno(2.2.4)$$
Put
$$v''_0:=u(s_1)v_0\q\h{with}\q u(s_1):=1/\bl(1+\al\1 s_1+\beta\1 h(s_1)\br).$$
Then $v'_0$ in the introduction can be obtained by modifying further this $v''_0$ as in (1.5). We have
$$v''_0=\ee_0+g_1\1\ee_1+\gt(s_1)\1\ee_2\mod\, V^{>1}\G_S,$$
with
$$\q g_1:=u(s_1)\bl(s_1+\gamma\1 h(s_1)\br),\q\gt:=u(s_1)h(s_1),$$
and this definition of $g_1$, $\gt$ is compatible with the one in the introduction. We see that $\gt$ is a function of $g_1$ (shrinking $S$ if necessary), since they are functions of $s_1$ with
$$g_1(0)=0,\q \dd g_1/\dd s_1(0)=1.$$
These imply that the ambiguity of $h$ in the expression of deformation is at most three dimensional (depending only on $\al$, $\beta$, $\gamma$).
\sk
We have for $j\ges 2$
$$pr_S^{(j)}(v''_0)=u(s_1)s_je_j=u(s_1)s_j\bl(\msum_{i=0}^{j-1}\,A_{j,i}\ee_i+\ee_j\br).$$
So the functions $g_j$ in the introduction are given for $j\ges 3$ by
$$g_j=u(s_1)s_j\q(3\les j\les r),
\leqno(2.2.5)$$
using (2.2.1) together with an argument in (1.5).
\sk
For the remaining case $j=2$, we have
$$pr_S^{(2)}(v''_0)=u(s_1)s_2\1e_2=u(s_1)s_2\bl(\beta\ee_0+\gamma\ee_1+\ee_2\br).$$
Using (2.2.2--3) for $k=1$ together with an argument in (1.5), we then get
$$g_2=\uu(s_1)s_2,
\leqno(2.2.6)$$
for some holomorphic function $\uu(s_1)$ of $s_1$ with $\uu(0)=1$.
\sk
These assertions imply that the $g_j$ ($j=1,\dots,r$) form a local coordinate system of $(S,0)$. Consequently no further condition is needed for the matrix $(A_{i,j})$. So the ambiguity of $h$ coming from automorphisms of $(S,0)$ and $\B_{h.S}$ coincides with the one coming from the choice of an opposite filtration $U$, and is exactly three dimensional. This finishes the proof of Theorem~1.
\msn
{\bf 2.3.~Remark.} It seems very difficult to extend the construction in this paper to the case the generation condition is satisfied, for instance, in the case $t\dt e_i=e_{i+1}$ ($i=0,\dots,r-1$) in the absolute case. There are also some cases where the generation condition is satisfied with monodromy semisimple, but the argument does not seem to work as well.
In the above first case, a versal deformation $\B_S$ contained in the constant deformation of Gauss-Manin system over $(S,0)=(\C^r,0)$ seems to be described as follows.
\sk
Define $N\in{\rm End}_\C(H)$ with $H:=\msum_{i=0}^r\,\C e_i$ by
$$Ne_i=e_{i+1}\,\,\,(i=0,\dots,r-1),\q Ne_r=0.$$
This can be extended to an endomorphism of $\G=\C\dtip[\dt]\otimes_{\C}H$ by the scalar extension under $\C\into\C\dtip[\dt]$, and the action of $t$ on $\G$ can be expressed by
$$(t\dt-i)=N\q\h{on}\q\dt^{-i}H\,\subset\,\G.
\leqno(2.3.1)$$
\sk
Put $(S,0):=(\C^r,0)$ with coordinates $s_1,\dots,s_r$. Set
$$\aligned v_j&:=\prod_{i=1}^r\exp(s_i\dt^{1-i}N^i)\dt^{-j}e_j\,\in\,\G_S:=\OO_S\dtip\otimes_{\C}H\q(j\in[0,r]),\\
\B_S&:=\bigoplus_{j=0}^r\OO_S\dtip v_j\,\,\subset\,\,\G_S.\endaligned
\leqno(2.3.2)$$
Here
$$[t,N]=[\dt^j,N]=[N,s_i]=[t,s_i]=[\dt^j,s_i]=0,\q[t,\dt^j]=-j\dt^{j-1}\,\,\,(j\in\Z).$$
We can verify that
$$\dti\dd_{s_j}v_0=\dt^{-j}N^j\prod_{i=1}^r\exp(s_i\dt^{1-i}N^i)e_0=\dt^{-j}N^jv_0=v_j\q(j\ges 1),
\leqno(2.3.3)$$
$$tv_0=t\,\prod_{i=1}^r\exp(s_i\dt^{1-i}N^i)e_0=\dti v_0+v_1+\sum_{i=2}^r(i-1)s_iv_i.
\leqno(2.3.4)$$
For $j=0,\dots,r$, we then get
$$\dti\dd_{s_i}v_j=\begin{cases}v_{i+j}&\h{if}\,\,\,i+j\les r,\\0&\h{if}\,\,\,i+j>r,\end{cases}
\leqno(2.3.5)$$
$$tv_j=(j+1)\dti v_j+v_{j+1}+\sum_{i=2}^{r-j}(i-1)s_iv_{i+j},
\leqno(2.3.6)$$
where $v_j:=0$ for $j>r$.
These may be shown, for instance, by using the following assertion:
\sk
If there are elements $A,B,C$ of a $\Q$-algebra such that
$$[A,B]=C,\q[B,C]=0,\q B^i=0\,\,\,(i\gg 0).$$
Then
$$[A,B^k]=k\1B^{k-1}C\,\,\,(k\in\N),\,\,\,[A,\exp B]=C\exp B.
\leqno(2.3.7)$$
\msn
{\bf 2.4.~Remark.} As a relative case, one may consider a one-parameter deformation of $\B$ in the introduction generated over $\OO_{S^{[1]}}\dtip$ by
$$\dt^{-j}\dd_{s_1}^jv^{[1]}_0\,\,\,(j=0,\dots,r)\q\h{with}\q
v^{[1]}_0:=e_0+s_1e_1+\msum_{i=2}^r\,h_i(s_1)e_i,
\leqno(2.4.1)$$
where $h_i(s_1)\in\C\{s_1\}s_1^i$ with $\dd_{s_1}^ih_i(0)\ne 0$ ($i\in\{2,\dots,r\}$), and $(S^{[1]},0):=(\C,0)$ with coordinate $s_1$. (Note that $N=0$ in the notation of (2.3).) In this case the generation condition together with the injectivity condition is satisfied, see \cite{HM}. This one-parameter family can be extended {\it canonically} to a family of Brieskorn lattices over $(S^{[r]},0)=(\C^r,0)$ contained in a constant family of Gauss-Manin systems and satisfying the condition
$$\dti\dd_{s_j}v^{[r]}_0\big|_{s=0}=\dt^{-j}\dd_{s_1}^jv^{[r]}_0\big|_{s=0}\q(j=2,\dots,r),
\leqno(2.4.2)$$
where $v^{[r]}_0$ is an extension of $v^{[1]}_0$ over $S^{[r]}$ defined as follows.
\sk
For $\nu=(\nu_2,\dots,\nu_r)\in\N^{r-1}$ with $\N:=\Z_{\ges 0}$, set
$$\sbo^{\nu}:=\mprod_{j=2}^r\,s_j^{\nu_j},\q\nu!:=\mprod_{j=2}^r\,\nu_j!,\q||\nu||=\msum_{j=2}^r\,j\1 \nu_j,\q|\nu|=\msum_{j=2}^r\,\nu_j.$$
Define $v^{[r]}_0$ by
$$v^{[r]}_0:=e_0+s_1e_1+\sum_{\nu\ges 0}\sum_{i=2}^rh_i^{(||\nu||)}(s_1)(\sbo^{\nu}/{\nu}!)\dt^{|\nu|-||\nu||}e_i,
\leqno(2.4.3)$$
where $h^{(k)}(s_1):=\dd_{s_1}^kh(s_1)$, and we have $\nu\ges\mu\iff\mu_j\ges\nu_j\,(\forall\,j\in\{2,\dots r\}$).
The Brieskorn lattice $\B_{S^{[r]}}$ on $S^{[r]}$ is generated over $\OO_{S^{[r]}}\dtip$ by
$$\dt^{-j}\dd_{s_1}^jv^{[r]}_0\q(j=0,\dots,r).$$
\sk
For $j=2,\dots,r$, we verify that
$$\aligned&\dt^{-j}\dd_{s_1}^jv^{[r]}_0=\sum_{\nu\ges 0}\sum_{i=2}^rh_i^{(||\nu||+j)}(s_1)(\sbo^{\nu}/{\nu}!)\dt^{|\nu|-||\nu||-j}e_i,\\
&\dti\dd_{s_j}v^{[r]}_0=\sum_{\nu\ges\ob_j}\sum_{i=2}^rh_i^{(||\nu||)}(s_1)(\sbo^{\,\nu-\ob_j}/(\nu-\ob_j)!)\dt^{|\nu|-||\nu||-1}e_i,\\
&\h{where}\q\q\ob_j=((\ob_j)_2,\dots,(\ob_j)_r)\in\N^{r-1}\q\h{with}\q(\ob_j)_i=\delta_{i,j}\,.\endaligned
\leqno(2.4.4)$$
Setting $\nu':=\nu+\ob_j$, we have
$$||\nu'||=||\nu||+j,\q|\nu'|-||\nu'||-1=|\nu|-||\nu||-j.$$
So we get the equalities
$$\dti\dd_{s_j}v^{[r]}_0=\dt^{-j}\dd_{s_1}^jv^{[r]}_0\q(j=2,\dots,r).
\leqno(2.4.5)$$
\sk
We have moreover
$$\aligned(t-\dti)v^{[r]}_0&=\sum_{\nu\ges 0}\sum_{i=2}^r\bl(||\nu||-|\nu|\br)h_i^{(||\nu||)}(s_1)(\sbo^{\nu}/{\nu}!)\dt^{|\nu|-||\nu||-1}e_i\\
&=\sum_{j=2}^r(j-1)s_j\dti\dd_{s_j}v^{[r]}_0.\endaligned
\leqno(2.4.6)$$
These imply that $\B_S^{[r]}$ is stable by the actions of $t$, $\dti\dd_{s_i}$ ($i=1,\dots,r$). Indeed, we have
$$\dti\dd_{s_i}(\dt^{-j}\dd_{s_1}^jv^{[r]})=\dt^{-j}\dd_{s_1}^j\dti\dd_{s_i}v^{[r]}=\dt^{-i-j}\dd_{s_1}^{i+j}v^{[r]},$$
and
$$\dt^{-j}\dd_{s_1}^jv^{[r]}\in V^{j+1}\G_S\subset\B_S\q\h{if}\q j>r,
\leqno(2.4.7)$$
where the last inclusion follows from the definition of $\B_S$. The argument is similar for the action of $t$ by using
$$(t-2\dti)\dti\dd_{s_i}v^{[r]}_0=\dti\dd_{s_i}(t-\dti)v^{[r]}_0.$$
\sk
In general we have
$$\dti\dd_{s_j}v^{[r]}_0\big|_{s=0}=\dt^{-j}\dd_{s_1}^jv^{[r]}_0\big|_{s=0}=
\sum_{i=2}^jh_i^{(j)}(0)\dt^{-j}e_i\,\ne\,h_j^{(j)}(0)\dt^{-j}e_j,
\leqno(2.4.8)$$
where the last inequality becomes an equality if the following condition is satisfied:
$$h_i(s_1)=a_is_1^i\q\h{with}\q a_i\in\C^*\,\,\,(i=2,\dots,r).
\leqno(2.4.9)$$
The above extension $v^{[r]}_0$ can be different from the one obtained by Theorem~(1.3) unless condition (2.4.9) is valid.
The uniqueness of extensions might be shown by assuming (2.4.9) although it does not seem necessarily easy to write it down carefully.
\msn
{\bf 2.5.~Remark.} It is possible to apply the argument in this paper to the case of {\it TE\1}-structures by using the partial Fourier transformation, that is, by replacing $\dti$, $t$ with $z_0$, $z_0^2\dd_{z_0}$ (up to a sign). Here $z_0$ is a coordinate of $\C$, and $\OO_S\dtip$ could be replaced by $\OO_{\C\times Z}$ with $Z$ a germ of a complex manifold so that one can consider finite free $\OO_{\C\times Z,0}$-modules $M$ with actions of $z_0^2\dd_{z_0}$, $z_0\dd_{z_i}$ ($i\ges 1$), where the $z_i$ ($i\ges 1$) are local coordinates of $Z$, see \cite{HM}, \cite{Sab}. It seems interesting to study a versal deformation of the Fourier transform $\B^F$ of $\B$ in the introduction, assuming it exists. Indeed, the arguments in this paper can be used to deduce a contradiction, assuming further $\dim Z=r+1$, see (3.8) below. Here the subset $Z_0$ of the base space $Z$ of the versal deformation on which $[z_0^2\dd_{z_0}]\in{\rm End}_{\OO_Z}(M/z_0M)$ is nilpotent is quite important. For the passage from deformations of Brieskorn lattices to the corresponding {\it TE\1}-structures, one can use {\it opposite filtrations} giving {\it algebraizations} of vector bundles, see {\it loc.~cit.} For the converse direction, we can simply take the scalar extension.
\bs\bs
\vbox{\centerline{\bf 3. Application}
\bsn
In this section we show the non-existence of a versal deformation of the Fourier transform $\B^F$ of $\B$ with dimension $r+1$ using the argument in the proof of Theorem~1.}
\msn
{\bf 3.1.~A deformation of $\B$.} For $\B$ in the introduction, there is a direct sum decomposition
$$\B=\mopl_{j=0}^r\,\B_j\q\,\,\,\h{with}\q\,\,\,\B_j=\C\dtip\1\dt^{-j}e_j.
\leqno(3.1.1)$$
Here each $\B_j$ is isomorphic to the Brieskorn lattice of $f_j:=\msum_{i=1}^{2j+2}\,x_i^2$ on $\C^{2j+2}$ ($j=0,\dots,r$). We have a 1-parameter deformation of $\B_j$ coming from the theory of Brieskorn lattices in the geometric case, that is, associated with a versal deformation of $f_j$. This is not contained in a trivial deformation of $\G$. Taking the direct sum of the pull-backs of these deformations under the projections from the product of their base spaces, we can get an $(r+1)$-parameter deformation of $\B$ satisfying the last property. This can be used to show that the subset $Z_0$ in Remark~(2.5) does not coincide with $Z$.
\msn
{\bf 3.2.~Period mapping into $\G$.} Let $\B$ be as above. For any deformation $\B_S$ of $\B$ on $S$ contained in the {\it constant\,\,} deformation $\G_S$, there is a unique $w_0\in\B_S$ corresponding to $e_0$ by Theorem~(1.3), where the opposite filtration $U$ is given by the $e_j$. This defines a {\it period mapping}
$$\sigma_S:S\to\G,
\leqno(3.2.1)$$
using the parallel translation induced by the triviality of the Gauss-Manin system $\G_S$ over $S$. This period mapping is quite different from the ones in the senses of K.~Saito and C.~Hertling. (In the latter, the Brieskorn lattice itself is used instead of the primitive form.) Here the target $\G$ is an infinite dimensional vector space, and this period mapping is essentially equivalent to considering {\it all} the coefficients of asymptotic expansions as in \cite{Va1}. By the construction of $w_0$ in the proof of Theorem~(1.3), however, the image of $\sigma_S$ is contained in a finite dimensional subspace
$$\C\,\dti\1e_0\oplus\bigoplus_{1\les j\les r,\,1\les k\les j}\C\,\dt^{-k}e_j\q\subset\q\G.
\leqno(3.2.2)$$
In the case $\B_S=\B_{h,S}$ for a holomorphic function $h$ of $s_1$, it is easy to see that $\sigma_S$ is {\it injective} and its image is a locally closed complex submanifold (shrinking $S$, $Z$ if necessary) by looking at the definition of the primitive form $v_0$ in the introduction.
\sk
The above period mapping is {\it compatible} with pull-backs by morphisms of base spaces $S$, that is, we have $\sigma_{S'}=\sigma_S\ssc g$ if $\sigma_{S'}$ is associated with the pull-back of the deformation by a morphism $g:S'\to S$.
\msn
{\bf 3.3.~Relation with a versal deformation of $\B^F$.} Let $\B$ and $\B_{h,S}$ be as above. Assume there is a versal deformation of the Fourier transform $\B^F$ of $\B$ with base space $(Z,0)$ and $\dim\,(Z,0)=r+1$ as in Remark~(2.5). Note that $\B^F$ is defined by using the basis $e_j$. The latter determines an opposite filtration $U$, and we there is a morphism of complex manifolds
$$\phi_h:(S,0)\to(Z,0),
\leqno(3.3.1)$$
such that $\B_{h,S}$ is isomorphic to the partial Fourier transform of the pull-back of the versal deformation of $\B^F$ by $\phi_h$. (Note that $\phi_h$ is not necessarily unique. We choose and fix one.)
\sk
Let $Z_0$ be the subset of $Z$ in Remark~(2.5) (that is, the residue of $z_0^2\dd_{z_0}$ is nilpotent on $Z_0$, where $1/z_0$ corresponds to $\dt$ by the Fourier transformation). We have the inclusion
$$\phi_h(S)\subset Z_0.
\leqno(3.3.2)$$
It follows from (3.1) that $Z_0\ne Z$. By the above hypothesis, this is equivalent to
$$\dim Z_0\les r.
\leqno(3.3.3)$$
On the other hand, we get by (3.2)
$$\dim\phi_h^{-1}(z)=0\q(\forall\,z\in Z_0).
\leqno(3.3.4)$$
(Indeed, if $\dim\phi_h^{-1}(z)>0$, then $\phi_h^{-1}(z)$ contains a locally closed smooth curve $C'$ on which the deformation $\B_{h,S}|_{C'}$ is {\it constant}, but this contradicts (3.2).) Since $\dim S=r$, these imply
$$\dim(Z_0,0)=r.
\leqno(3.3.5)$$
Moreover, shrinking $S$ and $Z$ if necessary, we may assume that $\phi_h$ is a {\it finite morphism} (by using the Weierstrass preparation theorem). Then $(\phi_h(S),0)$ is a {\it local irreducible component} of $(Z_0,0)$, and is a hypersurface. We have the following.
\msn
{\bf Proposition~3.4.} {\it With the above notation, $S$ is identified with the {\it normalization} of $\phi_h(S)$.}
\msn
{\it Proof.} It is enough to show that $\phi_h$ is {\it generically injective} (since $S$ is smooth and $\phi_h$ is a finite morphism). Let $U$ be a smooth non-empty Zariski-open subset of $\phi_h(S)$ over which $\phi_h$ is locally biholomorphic. Let $C$ be a locally irreducible curve on $\phi_h(S)$ passing through 0 and such that $C\setminus\{0\}\subset U$ (shrinking $S,Z$ if necessary).
\sk
Assume there is an irreducible component $D$ of $\phi_h^{-1}(C)$ such that the degree of the finite morphism $\phi_h|_D:D\to C$ is greater than 1. Let $\Dt$, $\Ct$ be the normalizations of $D$, $C$.  Consider the pull-back of $\B_{h,S}$ to $\Dt$ and also the partial Fourier transform of the pull-back to $\Ct$ of the versal deformation of $\B^F$. The former is isomorphic to the pull-back to $\Dt$ of the latter by a finite morphism of degree at least 2, and these are contained in the trivial deformation of Gauss-Manin systems on smooth curves $\Dt$, $\Ct$. We have the injectivity of the period mapping $\sigma_{\Dt}:\Dt\to\G$ by that of $\sigma_S:S\to\G$. However, $\sigma_{\Dt}$ is the composition of $\Dt\to\Ct$ with $\sigma_{\Ct}:\Ct\to\G$ by the compatibility of period mappings with pull-backs of base spaces as is explained at the end of (3.2). Hence it cannot be injective. So we get a contradiction.
\sk
If there are two curves $D_1$, $D_2$ over $C$, then we similarly get a contradiction by using the pull-back of the deformations to the normalizations of $C$, $D_1$, $D_2$ (since the period mapping is injective on $D_1\cup D_2\subset S$). Thus the generic injectivity is proved. This finishes the proof of Proposition~(3.4).
\msn
{\bf 3.5.~Observation.} In the notation and assumption of (3.3), assume there are holomorphic functions $h_1$, $h_2$ of $s_1$ together with $\phi_{h_1}$, $\phi_{h_2}$ such that
$$(\phi_{h_1}(S),0)=(\phi_{h_2}(S),0)\q\h{in}\q(Z_0,0).
\leqno(3.5.1)$$
Since $\phi_{h_i}$ is identified with the normalization ($i=1,2$) by Proposition~(3.4), there is a unique isomorphism
$$\rho:(S,0)\simto(S,0)\q\h{with}\q \phi_{h_1}=\phi_{h_2}\ssc\rho,
\leqno(3.5.2)$$
by the uniqueness of normalization. Since $\B_{h_i,S}$ is isomorphic to the partial Fourier transform of the pull-back of the versal deformation of $\B^F$ by $\phi_{h_i}$ ($i=1,2$), the last isomorphism implies an isomorphism
$$\B_{h_1,S}=\rho^*\B_{h_2,S}.
\leqno(3.5.3)$$
The argument in the introduction for the proof of Theorem~1 can be applied to this, and $h_1$ and $h_2$ are conjugate by the action of $\Gamma$ in the notation of (3.6) below.
\msn
{\bf 3.6.~Action of $\Gamma$.} Let $\Gamma$ be the group of lower triangular matrices $A$ of size $(r+1)\times(r+1)$ with $\C$-coefficients such that $A-{\rm Id}$ is nilpotent, where ${\rm Id}$ denotes the identity matrix. This group $\Gamma$ acts on the set $\Ob(F_0)$ of opposite filtrations $U$ to the Hodge filtration $F_0$ on $H$ such that
$$F_{0,p}H=\mopl_{i\les p}\,\C\,e_i\q(\forall\,p\in\Z).
\leqno(3.6.1)$$
This action is algebraic and simply transitive. As in (2.2) we have the action of $\Gamma$ on the deformations of Brieskorn lattices $\B_{h,S}$, which gives the action of $\Gamma$ on
$$\PS:=\bl\{h\in\C\{s_1\}\,\big|\,h(0)=h'(0)=0\br\}\,\bl(=\C\{s_1\}s_1^2\br).
\leqno(3.6.2)$$
This action is algebraic in the sense that it is naturally extended to the action on the formal power series $\PSh:=\C[[s_1]]s_1^2$, and moreover, if we have the expansions
$$h=\msum_{k\ges 2}\,c_ks_1^k,\q h^A=\msum_{k\ges 2}\,c^A_ks_1^k,
\leqno(3.6.3)$$
then, setting as in (2.2)
$$\al:=A_{1,0},\q\beta:=A_{2,0},\q\gamma:=A_{2,1},
\leqno(3.6.4)$$
we have
$$c^A_k\in\C[\al,\beta,\gamma,c_2,\dots,c_k]\q(\forall\,k\ges 2).
\leqno(3.6.5)$$
Here the action of $A\in\Gamma$ on $h$ is denoted by $h^A$, since the action of $\Gamma$ on $\PS$ as is defined in (2.2) is contravariant, see (3.7) below.
\sk
For instance, the coefficients $b_k$ of the inverse function $x=\msum_{k\ges 1}\,b_k\1 y^k\in\C\{y\}$ of a convergent power series $y=\msum_{k\ges 1}\,a_kx^k\in\C\{x\}$ with $a_1=1$ can be obtained as a polynomial of $a_2,\dots,a_k$ by increasing induction on $k\ges 1$ using the equalities
$$b_k\1x^k+\msum_{i=1}^{k-1}\,b_i\,\bl(\msum_{j=1}^k\,a_jx^j\br)^i=x\mod(x^{k+1}).
\leqno(3.6.6)$$
\sk
So we get the induced algebraic action of $\Gamma$ on affine spaces
$$\Ab_k:=\C\{s_1\}s_1^2/\C\{s_1\}s_1^{k+2}=\C^k\q(\forall\,k\ges 1),
\leqno(3.6.7)$$
in a compatible way with the natural projections
$$\PS\onto\Ab_k.
\leqno(3.6.8)$$
These are compatible with the isomorphism
$$\PSh=\rlap{\raise-10pt\h{$\,\,\,\scriptstyle k$}}\rlap{\raise-6pt\h{$\,\leftarrow$}}{\rm lim}\,\Ab_k.
\leqno(3.6.9)$$
\msn
{\bf 3.7.~Contravariant action.} We have a contravariant action of the group $\Gamma$ on $\PS$ in the notation of (3.6). We can verify this directly by considering the action of $\Gamma$ on the set of pairs $(g_1,\gt)$ with $g_1\in\C\{s_1\}s_1$, $\gt\in\C\{s_1\}s_1^2$, and $\dd g_1/\dd s_1(0)=1$ as in (2.2). Indeed, in the notation of (2.2), setting
$$Ae_i=\msum_j\,A_{i,j}e_j,
\leqno(3.7.1)$$
we get the action of $A=(A_{i,j})\in\Gamma$ on
$$e_0+s_1\1 e_1+h(s_1)\1 e_2
\leqno(3.7.2)$$
by 
$$\bl(1+\al\1 s_1+\beta\1 h(s_1)\br)\1 e_0+\bl(s_1+\gamma\1 h(s_1)\br)\1 e_1+h(s_1)\1 e_2,
\leqno(3.7.3)$$
where $\al$, $\beta$, $\gamma$ are as in (3.6) (or (2.2)).
Dividing it by $1+\al\1 s_1+\beta\1 h(s_1)$, we get as in (2.2)
$$e_0+g_1\1 e_1+\gt\1 e_2,
\leqno(3.7.4)$$
with
$$\q g_1:=\frac{s_1+\gamma\1 h(s_1)}{1+\al\1 s_1+\beta\1 h(s_1)},\q
\gt:=\frac{h(s_1)}{1+\al\1 s_1+\beta\1 h(s_1)}.
\leqno(3.7.5)$$
Applying further the action of $A'$ on this, we get the following:
$$\aligned&\biggl(1+\frac{\al'\1s_1+\al'\1\gamma\1 h(s_1)}{1+\al\1 s_1+\beta\1 h(s_1)}+\frac{\beta'\1h(s_1)}{1+\al\1 s_1+\beta\1 h(s_1)}\biggr)\1 e_0\\
&+\biggl(\frac{s_1+\gamma\1 h(s_1)}{1+\al\1 s_1+\beta\1 h(s_1)}
+\frac{\gamma'\1h(s_1)}{1+\al\1 s_1+\beta\1 h(s_1)}\biggr)\1 e_1\\
&+\biggl(\frac{h(s_1)}{1+\al\1 s_1+\beta\1 h(s_1)}\biggr)\1 e_2\endaligned
\leqno(3.7.6)$$
Multiplied by $1+\al\1 s_1+\beta\1 h(s_1)$, this becomes
$$\aligned&\bl(1+(\al+\al')s_1+(\beta+\beta'+\al'\1\gamma)h(s_1)\br)\1 e_0\\
&+\bl(s_1+(\gamma+\gamma')h(s_1)\br)\1 e_1\\
&+h(s_1)\1 e_2\endaligned.
\leqno(3.7.7)$$
On the other hand, we have
$$\begin{pmatrix}1&0&0\\\al&1&0\\\beta&\gamma&1\end{pmatrix}
\begin{pmatrix}1&0&0\\\al'&1&0\\\beta'&\gamma'&1\end{pmatrix}
=\begin{pmatrix}1&0&0\\\al+\al'&1&0\\\beta+\beta'+\gamma\al'&\gamma+\gamma'&1\end{pmatrix}
\leqno(3.7.8)$$
So we get a contravariant action. Note that this property is closely related to the following assertion:
$$\h{If $\,Ae_i=\msum_j\,A_{i,j}e_j\,$ and $\,A'e_j=\msum_k\,A'_{j,k}e_k,\,$ then $\,A'Ae_i=\msum_{j,k}\,A_{i,j}A'_{j,k}e_k$.}
\leqno(3.7.9)$$
\msn
{\bf 3.8.~Non-existence of a versal deformation of $\B^F$ with dimension $r+1$.} In the notation and assumption of (3.3) and (3.6), let $\PS/\Gamma$ be the quotient of $\PS$ by the action of $\Gamma$, and ${\rm Irr}(Z,0)_r$ be the set of irreducible components of $(Z_0,0)$ with dimension $r$. We denote by ${\rm Irr}(Z,0)'_r$ the subset consisting of $(\phi_h(S),0)\in{\rm Irr}(Z,0)_r$ for some $\phi_h$. (Note that $\phi_h$ is not necessarily uniquely determined by $h$ as is remarked after (3.3.1).) The image of $h$ in $\PS/\Gamma$ is denoted by $[h]$. (It may be viewed as an {\it invariant} of the isomorphism class of the deformation $\B_{h,S}$.)
By the observation in (3.5), there is a well-defined surjective map
$${\rm Irr}(Z,0)'_r\ni(\phi_h(S),0)\mapsto[h]\in\PS/\Gamma.
\leqno(3.8.1)$$
This is, however, a contradiction, since ${\rm Irr}(Z,0)'_r$ is finite. So the versal deformation as in (3.3) cannot exist.

\end{document}